\documentclass[preprint,amscd,amsmath,amssymb,verbatim]{revtex4}

\usepackage{graphicx}
\usepackage{textcomp}
\usepackage[OT1,T1]{fontenc}

\begin{document}

\title{Riemann Hypothesis: The Riesz-Hardy-Littlewood wave in the
long wavelength region}
\author{Stefano Beltraminelli}
\email{stefano.beltraminelli@ti.ch}
\affiliation{CERFIM, Research Center for Mathematics and Physics,
PO Box 1132, 6600 Locarno, Switzerland}
\affiliation{ISSI, Institute for Scientific and Interdisciplinary
Studies, 6600 Locarno, Switzerland}\author{Danilo Merlini}
\email{merlini@cerfim.ch}
\affiliation{CERFIM, Research Center for Mathematics and Physics,
PO Box 1132, 6600 Locarno, Switzerland}
\affiliation{ISSI, Institute for Scientific and Interdisciplinary
Studies, 6600 Locarno, Switzerland}\label{I1}\label{I2}
\date{\today}
\begin{abstract}
We present the results of numerical experiments in connection with
the Riesz and Hardy-Littlewood criteria for the truth of the Riemann
Hypothesis (RH). The coefficients $c_{k}$ of the Pochammer's expansion
for the reciprocal of the Riemann Zeta function, as well as the
``critical functions'' $c_{k}k^{a}$ (where {\itshape a} is some
constant), are analyzed at relatively large values of {\itshape
k}. It appears an oscillatory behaviour (Riesz-Hardy-Littlewood
wave). The amplitudes and the wavelength of the wave are compared
with an analytical treatment concerning the wave in the asymptotic
region. The agreement is satisfactory. We then find numerically
that in the large $\beta$ limit too, the amplitudes of the waves
appear to be bounded. For a special case the numerical experiments
are performed up to larger values of {\itshape k}, i.e $k={10}^{9}$
and more. The analysis suggests that RH may barely be true and an
absolute bound for the amplitudes of the waves in all cases\ \ should
be given by $|\frac{1}{\zeta ( \frac{1}{2}+\epsilon ) }-1|$, with
$\epsilon$ arbitrarily small positive, i.e. equal to 1.68\ldots
.
\end{abstract}
\maketitle

\section{Introduction}

Following recent works concerning the study of some well known functions
appearing in the original criteria of Riesz, Hardy and Littlewood
for the possible truth of the Riemann Hypothesis (RH), there is
new interest in the direction of numerical experiments, where the
calculations use the ideas of some recent works on the subject.
These concern the expansion of the reciprocal of the Riemann Zeta
function in terms of the so called Pochammer's polynomials $P_{k}$,
whose coefficients $c_{k}$ play a central role also in the asymptotic
region of very large {\itshape k }\cite{1,2,3,4}. For new zero free
regions of the Zeta function, in the context of a rigorous treatment
with the M\"untz formula, the reader may consult a recent work by
Albeverio and Cebulla \cite{5}.

Here we are concerned with the discrete version of the Riesz criterion
which has also been studied numerically: the first numerical experiments
for values of {\itshape k} up to 100'000 have been announced and
reported for the Riesz case in \cite{2,6}. It has been found that
the function $c_{k}$ has a oscillatory behaviour in a region of
relatively high {\itshape k}'s, in agreement with an asymptotic
formula given by Baez-Duarte \cite{2}. The agreement appears satisfactory
even if only the contribution of the first non trivial zero of the
Zeta function located at $s=\frac{1}{2}+i 14.134725$ in the complex
plane has been used.

In a previous work \cite{7} a two parameter family (parameter $\alpha$
and $\beta$) of Pochammer's polynomials was introduced. This allowed
the starting investigation of $c_{k}$ at low values of {\itshape
k}, but in various cases and in the so called ``strong coupling''
regime (high $\beta$). After the initial study at low {\itshape
k}, our computations using the formula containing the M\"obius function
were easily extended to larger and larger {\itshape k} (up to a
billion) in the strong coupling limit, with the appearance of macro-oscillations
in $c_{k}$ extending to larger {\itshape k}. This is a symptom that
using such a limit the RH may eventually barely be true\cite{7}.
In this work we continue the numerical experiments also using our
Poisson formula established in \cite{7} which is well suited for
numerical purposes.

After the formulation of the model in Section 2, we then compute
in Section 3 the amplitudes of what we call the Riesz-Hardy-Littlewood
wave, which contains arbitrarily scales in few of the two parameters
$\alpha$ and $\beta$ at our disposal. Using the Baez-Duarte formula,
we then present our results for different models up to values of
{\itshape k} equal to one billion and observe oscillations in all
cases (Section 4). The agreement with the asymptotic formula of
Baez-Duarte is satisfactory. Then, in Section 5, we concentrate
the study in more details by considering a special new model already
proposed in \cite{7} where $\alpha =\frac{7}{2}$ and $\beta$ is
increasing starting with the value equal to 4. The results show
in a concrete way the ``transition'' from the low coupling to the
``strong coupling regime'': at low values of $\beta$ ($\beta =4$)
we obtain up to 7 oscillation with values of {\itshape k} extending
up to a billion. These start to deform continuously with increasing
values of $\beta$ approaching the infinite $\beta$ limit. In such
a regime, the wave is absorbed in a macroscopic region with an amplitude
whose strength should be finite as already noted in \cite{7}.

In the context of validity of our numerical results, our analysis
gives further indication that the RH may barely be true as indicated
by our two parameter models in the week as well as in the ``strong
coupling regime'' (Section 5). Moreover, the possibility that in
an ideal numerical experiment (using an arbitrarily large but finite
maximum value of {\itshape n}, say {\itshape N} in the formula with
the M\"obius function) the amplitude of the waves at finite $\beta$
values should decrease, is commented in Appendix.

\section{The model}\label{XRef-Section-34105336}

Following recent treatments\cite{1,2,7}, a possible expansion of
the reciprocal of the Zeta function i.e. ${\zeta ( s) }^{-1}$ in
terms of the so called Pochammer's polynomials $P_{k}$, with two
parameters $\alpha$ ($\alpha >1$) and $\beta$ is this one:
\begin{equation}
\frac{1}{\zeta ( s) }=\sum \limits_{k=0}^{\infty }c_{k}( \alpha
,\beta ) P_{k}( s,\alpha ,\beta ) 
\end{equation}

\noindent where
\begin{align}
P_{k}( s,\alpha ,\beta ) &=\prod \limits_{r=1}^{k}\left( 1-\frac{\frac{s-\alpha
}{\beta }+1}{r}\right) 
\\%
c_{k}( \alpha ,\beta ) &=\sum \limits_{n=1}^{\infty }\frac{\mu (
n) }{n^{\alpha }}{\left( 1-\frac{1}{n^{\beta }}\right) }^{k}
\end{align}

\noindent and $P_{k}( 0,\alpha ,\beta ) =1$. 

In (3) the M\"obius function of argument {\itshape n} is given by:
\[
\mu ( n) =\begin{cases}
1, & \mathrm{if} \ n=1 \\
{\left( -1\right) }^{k}, & \mathrm{if} \ n \ \mathrm{is} \ \mathrm{a}
\ \mathrm{product} \ \mathrm{of} \ k \ \mathrm{distinct} \ \mathrm{primes}
\\
0, & \mathrm{if} \ n \ \mathrm{contains} \ \mathrm{a} \ \mathrm{square}
\\
\end{cases}
\]

If $s=\rho +i t $ is a complex variable and argument of the Riemann
Zeta function one has for $\mathfrak{R}( s) =\rho >1$:
\begin{equation}
\frac{1}{\zeta ( s) }=\sum \limits_{n=1}^{\infty }\frac{\mu ( n)
}{n^{s}}
\end{equation}

Another explicit formula for the $c_{k}( \alpha ,\beta ) $ is obtained
from (3) using the binomial coefficients and reads: 
\begin{equation}
c_{k}( \alpha ,\beta ) =\sum \limits_{j=0}^{k}{\left( -1\right)
}^{j}\binom{k}{j}\frac{1}{\zeta ( \alpha +\beta  j) }
\end{equation}

As $\beta$ is increasing, one may also use (especially) in the context
of numerical experiments, the formula recently obtained \cite{7}
and given by:
\begin{equation}
c_{k}( \alpha ,\beta ) =\sum \limits_{n=1}^{\infty }\frac{\mu (
n) }{n^{\alpha }}e^{-\frac{k}{n^{\beta }}}
\end{equation}

In such an approximation we have that
\begin{equation}
c_{k}( \alpha ,\beta ) =\sum \limits_{p=0}^{\infty }c_{p}( \alpha
,\beta ) \frac{k^{p}}{p!}e^{-k}
\end{equation}

\noindent which shows the emergence of a Poisson like distribution
for the coefficients $c_{k}( \alpha ,\beta ) $. This should be a
very satisfactory approximation in the limit of relative large values
of $\beta$ \cite{7}. We recall that an important inequality due
to Baez-Duarte \cite{2}, concerning the Pochammer's polynomials
of complex argument {\itshape z} is given by:
\begin{equation}
\left| P_{k}( z) \right| \leq C k^{-\mathfrak{R}( z) }
\end{equation}

The above inequality applied to our two parameter family of Pochammer's
polynomials with complex argument $z=\frac{s-\alpha }{\beta }+1$
gives:
\begin{equation}
\left| P_{k}( s,\alpha ,\beta ) \right| \leq C k^{-\frac{\rho -\alpha
}{\beta }+1}
\end{equation}

So that $\zeta ( s) $ in (1) will be different from zero and thus
the RH will be true for $\mathfrak{R}( s) =\rho >\frac{1}{2}$ if
the $c_{k}$ decay, at large {\itshape k}, as (see \cite{7}): 
\begin{equation}
\left| c_{k}\right| \leq A k^{-\frac{\alpha -\rho }{\beta }}
\end{equation}

We will also consider the ``critical function''\,
\begin{equation}
\psi ( k;\alpha ,\beta ,\rho ) :=c_{k}k^{\frac{\alpha -\rho }{\beta
}}
\end{equation}

\noindent which from (10) is expected to be bounded by a constant
{\itshape A}.

We now recall two original cases given in pionnering works by Riesz
\cite{8} and by Hardy-Littlewood \cite{9}. Setting $\rho =\frac{1}{2}$
in (10), for $\alpha =\beta =2$ (Riesz case) we have that $|c_{k}|\leq
A k^{-\frac{3}{4}}$ and for $\alpha =1$, $\beta =2$ (Hardy-Littlewood
case) $|c_{k}|\leq A k^{-\frac{1}{4}}$. Other interesting cases
for which we will carry out intensive numerical experiments to be
presented below are summarized in the Table I.
\begin{table}
\caption{The expected decay of $c_{k}$ for different values of $\alpha$
and $\beta$} 
\begin{ruledtabular}

\begin{tabular}{ccccl}
$\alpha$ & $\beta$ & $\rho$ & decay of $|c_{k}|$ & Note\\
\hline
2 & 2 & $\frac{1}{2}$ & $k^{-\frac{3}{4}}$  & The case of Riesz\\
1 & 2 & $\frac{1}{2}$ & $k^{-\frac{1}{4}}$ & The case of Hardy-Littlewood\\
2 & 6 & $\frac{1}{2}$ & $k^{-\frac{1}{4}}$ & Same decay as the Hardy-Littlewood
case but numerically more convenient\\
$\frac{7}{2}$ & 4 & $\frac{1}{2}$ & $k^{-\frac{3}{4}}$ & Same decay
as the Riesz case, intensive calculations are given below\\
3 & 3 & $\frac{1}{2}$ & $k^{-\frac{5}{6}}$ & If the Zeta function
has no zero for $\rho >\frac{3}{4}$ then $c_{k}( 3,3) $

should decay at least as $k^{-\frac{3}{4}}$\ \ \\
4 & 4 & $\frac{1}{2}$ & $k^{-\frac{7}{8}}$ & Since from the Prime
number theorem there is no zero for $\rho >1$ \\ 
 \ &  \ &  \ &  \ & the $c_{k}( 4,4) $ decays at least as $k^{-\frac{3}{4}}$\\
2 & 4 & $\frac{1}{2}$ & $k^{-\frac{3}{8}}$ & Another interesting
case for calculations
\end{tabular}
\end{ruledtabular}
\end{table}

A limiting delicate case analyzed in \cite{7} is the one where $\alpha
=\frac{1}{2}+\delta $ and $\beta$ grows to infinity. Here of course
we do not have absolute convergence to ${\zeta ( s) }^{-1}$ ($c_{k}$
may nevertheless be analyzed) and from (10) we have that the $c_{k}$
should be smaller then a constant for all {\itshape k}. This is
what we verified with numerical experiments (not presented here)
with values of {\itshape k} up to a billion. The value of the constant
has been proposed in our previous work \cite{7} and the conjecture
was that $|c_{k}|\leq |\frac{1}{\zeta ( \frac{1}{2}) }-1|\cong 1.68477$.
However the situation is delicate ($\alpha <1$) since Littlewood
\cite{10} has shown that, assuming RH is true, $\sum \limits_{n=1}^{\infty
}\frac{\mu ( n) }{n^{\frac{1}{2}+\epsilon }}$ is convergent for
all $\epsilon$ strictly greater than zero.

The general situation is that the ``critical function'' $c_{k}k^{\frac{\alpha
-\rho }{\beta }}$ should be bounded by a constant in absolute value
as $k\overset{ }{\rightarrow }\infty $. In fact the function starts
at zero for $k=0$, reaches a minimum, then starts to increase and
then begins to oscillate with a ``constant `` amplitude as $k\overset{
}{\rightarrow }\infty $ as we will see in the experiments. In a
previous work \cite{7} we have analyzed $c_{k}$ in various cases
but only for moderately values of {\itshape k}, i.e for {\itshape
k} not exceeding 1000, with exception of some cases at large values
of $\beta$, where {\itshape k} reached the value of a half billion.
$c_{k}$ was found to have only negative values in the range considered
and increasing with {\itshape k}. Presently we know of recent numerical
experiments \cite{6} in the Riesz case carried out by Maslanka ({\itshape
k} up to 100'000) and Wolf ({\itshape k} up to 200'000). These calculations
show that $c_{k}$ become of oscillatory type, thus assuming positive
and negative values with an amplitude which appears constant in
the range considered. In fact two or three oscillations with a wavelength
related in first approximation to the first zero of the Riemann
Zeta function may be seen. Here it should be remarked that this
situation for the Riesz case is not in contraddiction with our strong
coupling limit ($\beta$ large) cited above (see discussion below
for the case $\alpha =\frac{7}{2}$ and $\beta$ increasing). 

\ In few of these new finding, we want first analyze (in an analytical
context) such a behaviour and we call this general phenomena the
Riesz-Hardy-Littlewood wave. This will be analyzed using an interesting
result of Baez-Duarte, i.e. an expression giving $c_{k}$ for $k\overset{
}{\rightarrow }\infty $.

\section{The Riesz-Hardy-Littlewood wave}

 For the Riesz case, in connection with the Mellin inversion formula,
the Riesz function is given (see \cite{8} and \cite{11}) explicitly
by:
\begin{equation}
F( x) =\sum \limits_{k=1}^{\infty }\frac{{\left( -1\right) }^{k+1}x^{k}}{\left(
k-1\right) !\zeta ( 2k) }
\end{equation}

Using the calculus of residues {\itshape F(x)} is obtained by an
integration and is given by:
\begin{equation}
F( x) =\frac{i}{2\pi }\operatorname*{\int }\limits_{a-i \infty }^{a+i
\infty }\frac{\Gamma ( 1-s) x^{s}}{\zeta ( 2s) }ds\ \ \ \ 
\end{equation}

\noindent where $\frac{1}{2}<a<1$.

Now, recently Baez-Duarte \cite{2}, with an ingenious method found
in particular an expression for the reciprocal of the Pochammer
polynomial given by:
\begin{equation}
\frac{1}{P_{k}( s) }=\sum \limits_{j=1}^{k}{\left( -1\right) }^{j}\binom{k}{j}\frac{j}{s-j}
\end{equation}

\noindent where uniformely on compact subsets one has:\ \ \ \ \ 
\begin{equation}
\operatorname*{\lim }\limits_{k\,\rightarrow \:\infty }P_{k}( s)
k^{s}=\frac{1}{\Gamma ( 1-s) }
\end{equation}

\noindent and he was able to obtain an explicit formula connecting
$c_{k}$ and the set of all trivial and non trivial zeros ({\itshape
z} denote the complex Zeta zeros, $z=\frac{1}{2}+i t$) under the
assumption of simple zeros. For the Riesz case the expression is
given by:
\begin{equation}
-2k c_{k-1}=\sum \limits_{\mathfrak{I}( z) }^{ }\frac{1}{\zeta ^{\prime
}( z) P_{k}( \frac{z}{2}) }
\end{equation}

\noindent for sufficiently large {\itshape k}. It should be said
that formula (16) of Baez-Duarte is very nice and may be used to
control our numerical computations at large {\itshape k} to be presented
below. Apparently (16), with some precautions, may be extended to
the general case with parameters $\alpha$, $\beta$ and should read:
\begin{equation}
-\beta  k c_{k-1}=\sum \limits_{\mathfrak{I}( z) }^{ }\frac{1}{\zeta
^{\prime }( z) P_{k}( \frac{s-\alpha }{\beta }+1) }
\end{equation}

To obtain an asymptotic value for the amplitude of the Riesz-Hardy-Littlewood
wave, we use (15) in (17):
\begin{equation}
-\beta  k c_{k-1}=\sum \limits_{\mathfrak{I}( z) }^{ }\frac{k^{\frac{i
t}{\beta }}k^{\frac{\frac{1}{2}-\alpha }{\beta }}}{\zeta ^{\prime
}( z) }
\end{equation}

\noindent where $t=\mathfrak{I}( z) $. In the limit of large {\itshape
k}, one may neglect the contribution of the trivial zeros \cite{2}.
For the ``critical function'' we then have the following expression:
\begin{equation}
k^{\frac{\alpha -\frac{1}{2}}{\beta }}c_{k}\cong -\frac{1}{\beta
}\sum \limits_{\mathfrak{I}( z) }^{ }\frac{k^{\frac{i t}{\beta }}\Gamma
( -\frac{\frac{1}{2}+i t-\alpha }{\beta }) }{\zeta ^{\prime }( z)
}=:\overline{\psi }( k;\alpha ,\beta ,\frac{1}{2}) 
\end{equation}

\noindent for large {\itshape k}. To prepare the comparison of (19)
with the numerical results we write explicitly (19) for the various
cases we will treat. In order to obtain an estimate for the amplitude
of the wave in the long wavelength limit ({\itshape k} large) we
will use here only the first zero of the Riemann Zeta function up
to 10 decimals ($\beta$ small).
\begin{table}
\caption{The ``amplitude'' of $\psi$ for different values of $\alpha$
and $\beta$} 
\begin{ruledtabular}

\begin{tabular}{lccc}
$\alpha$ & $\beta$ & The function & The amplitude\\
\hline
2 & 2 & $|k^{\frac{3}{4}}c_{k}|=|\psi ( k;2,2,\frac{1}{2}) |$ &
0.000078\\
1 & 2 & $|k^{\frac{1}{4}}c_{k}|=|\psi ( k;1,2,\frac{1}{2}) |$ &
0.0000292558\\
2 & 6 & $\text{$|k^{\frac{1}{4}}c_{k}|=|\psi ( k;2,6,\frac{1}{2})
|$}$ & 0.0210433\\
$\frac{7}{2}$ & 4 & $|k^{\frac{3}{4}}c_{k}|=|\psi ( k;\frac{7}{2},4,\frac{1}{2})
|$ & 0.008411\\
3 & 3 & $|k^{\frac{5}{6}}c_{k}|=|\psi ( k;3,3,\frac{1}{2}) |$ &
0.0021562\\
4 & 4 & $|k^{\frac{7}{8}}c_{k}|=|\psi ( k;4,4,\frac{1}{2}) |$ &
0.00984936\\
2 & 4 & $|k^{\frac{3}{8}}c_{k}|=|\psi ( k;2,4,\frac{1}{2}) |$ &
0.0052445
\end{tabular}
\end{ruledtabular}
\end{table}

\ \ The upper bounds for the amplitude of the waves above, will
be compared with the results of the numerical experiments performed
for the various cases using (3).

\section{Numerical experiments}

We now present the results of our numerical experiments which was
carried out in more cases using the M\"obius function in (3) up
to $n={10}^{6}$. We calculated $c_{k}$ until $k={10}^{6}$ or $k={10}^{9}$
with a scaling factor of 2500 for the {\itshape k}-axis. These will
be compared with the upper bound for the amplitude of the waves
of Section 3. The general situation is that for moderately values
of {\itshape k }(until some tausend) the wave given by the experimental
results start with zero amplitude, after a minimum with a negative
value, increases and seems to stabilize at large values of {\itshape
k} with oscillations displaced at larger and larger wavelength (proportional
to $\log ( k) $) and with an amplitude which seems to saturate to
a constant value (given in a good approximation) by the upper bound
(19). Below (Figure 1) we first give the plots of the wave for the
Riesz case ($\alpha =\beta =2$). As remarked in \cite{2}, the first
intensive calculations with very high precision up to $k=100'000$
(by K. Maslanka) and up to $k=200'000$ (by M. Wolf) indicated the
appearance of oscillations with the first one in the region $k=20'000$.
Our results obtained with (3) confirm for such values the asymptotic
limit for the wave with an amplitude in agreement with the bound
obtained above ($A\cong 0.000078$). Notice that the minimun for
low {\itshape k} values is 0.4 in absolute value as found in a previous
work, is much bigger then {\itshape A,} which concerns only the
asymptotic region of the wave, thus no disagreement!
\begin{figure}[h]
\begin{center}
\includegraphics{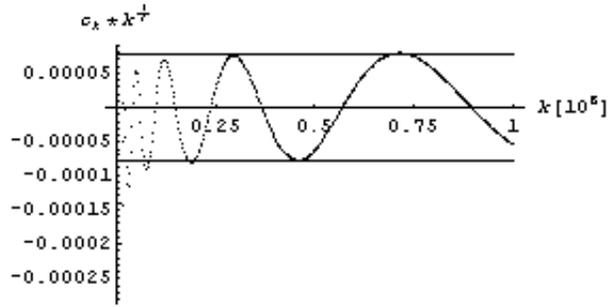}

\end{center}
\caption{The wave $k^{\frac{3}{4}}c_{k}$ for the Riesz case $\alpha
=\beta =2$}

\end{figure}
\begin{figure}[h]
\begin{center}
\includegraphics{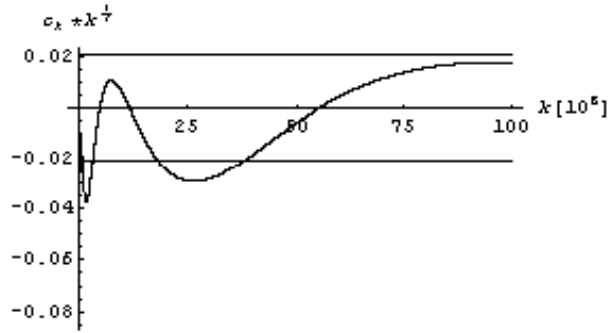}

\end{center}
\caption{The wave $k^{\frac{1}{4}}c_{k}$ for the case $\alpha =2,\beta
=6$}

\end{figure}

Figure 2 and Figure 3 concern two cases of special interest since
the decays are expected to be the same as for the Hardy-Littlewood
case and for the Riesz case. In both cases there is agreement with
the bound $A\cong 0.0210433$ and $A\cong 0.008411$ given above but
the amplitudes are respectively 1000 and 100 time bigger than in
the former cases.
\begin{figure}[h]
\begin{center}
\includegraphics{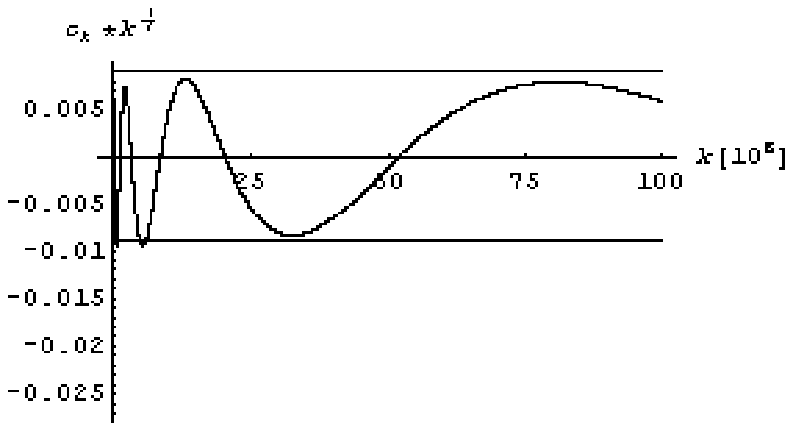}

\end{center}
\caption{The wave $k^{\frac{3}{4}}c_{k}$ for the case $\alpha =\frac{7}{2},\beta
=4$}

\end{figure}

In Figure 4 and Figure 5 we give the plots of $k^{\frac{5}{6}}c_{k}(
3,3) $ and $k^{\frac{7}{8}}c_{k}( 4,4) $ where the amplitudes are
found to be in agreement with the theoretical upper bounds given
above in Table 2, too.
\begin{figure}[h]
\begin{center}
\includegraphics{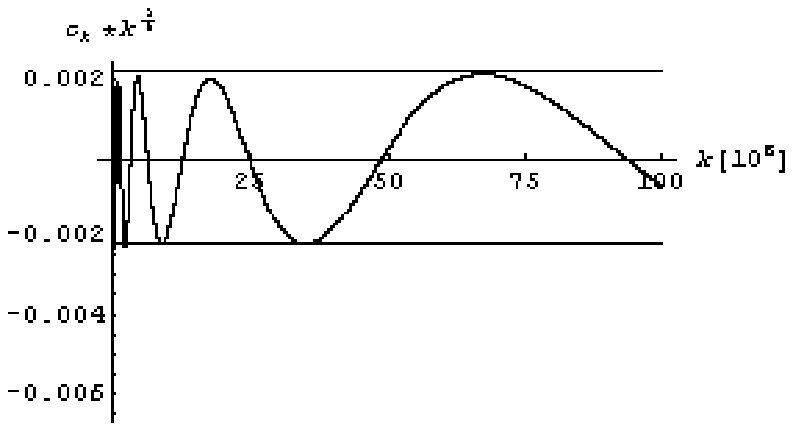}

\end{center}
\caption{The wave $k^{\frac{5}{6}}c_{k}$ for the case $\alpha =\beta
=3$}

\end{figure}
\begin{figure}[h]
\begin{center}
\includegraphics{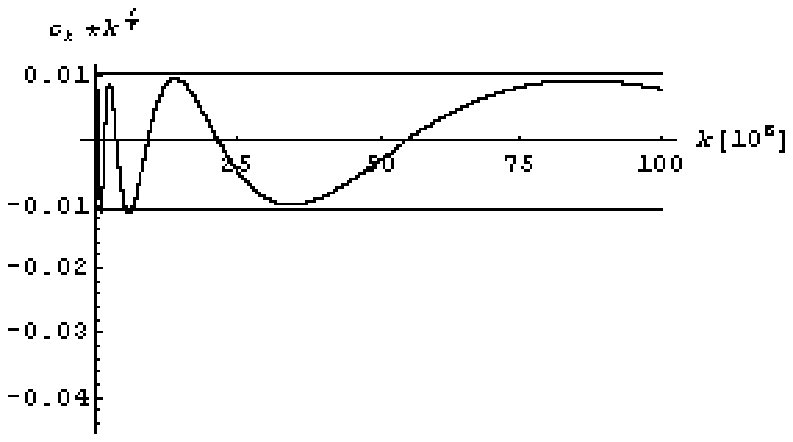}

\end{center}
\caption{The wave $k^{\frac{7}{8}}c_{k}$ for the case $\alpha =\beta
=4$}

\end{figure}

The next special case is the one with $\alpha =2$ and $\beta =4$.
Again, the experimentally detected amplitude agrees well with the
theoretical bound given above, i.e $A=0.0052445$.
\begin{figure}[h]
\begin{center}
\includegraphics{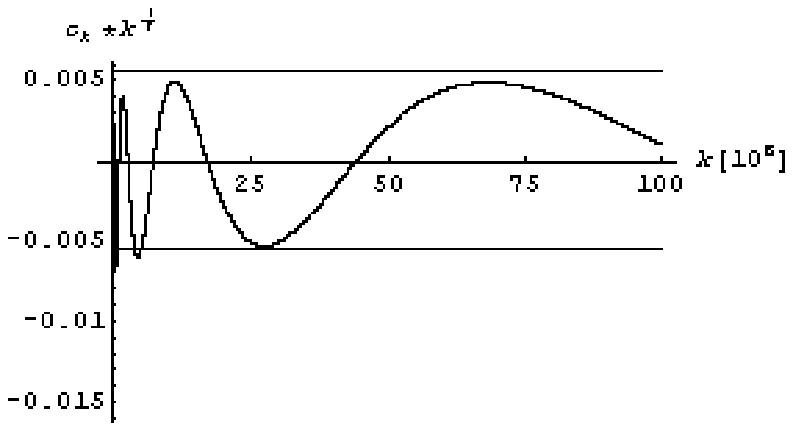}

\end{center}
\caption{The wave $k^{\frac{3}{8}}c_{k}$ for the case $\alpha =2,\beta
=4$}

\end{figure}

As a further illustration we compare the wave $k^{\frac{3}{8}}c_{k}$
with the asymptotic approximation wave given by the wave $\psi (
k;2,4,\frac{1}{2}) $ (case $\alpha =2,\beta =4$). In the range for
{\itshape k} from $1\cdot {10}^{6}$ to $10\cdot {10}^{6}$ the two
waves appear to be walking close together arm in arm (Figure 7).
Notice that in the approximation we considered only the contribution
of the first zero given by $z=\frac{1}{2}+i 14.134725141$ which
appears dominant for low values of $\beta$.
\begin{figure}[h]
\begin{center}
\includegraphics{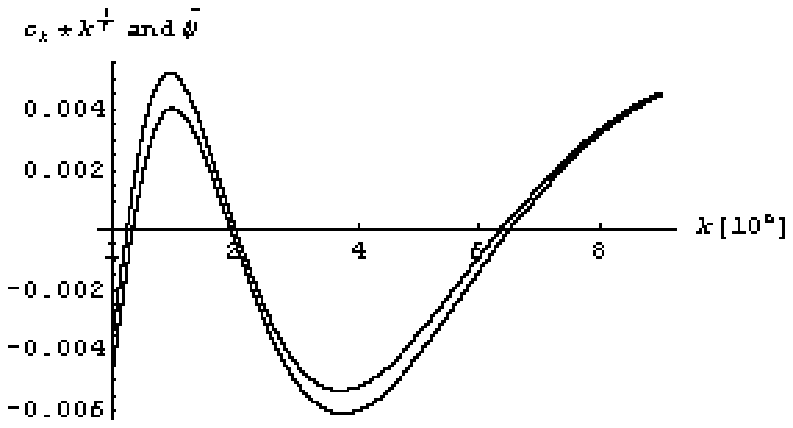}

\end{center}
\caption{The wave $k^{\frac{3}{8}}c_{k}$ (lowest curve) and the
approximation $\overline{\psi }$ (highest curve)}

\end{figure}

\section{The case $\alpha =\frac{7}{2}$ and $\beta$ increasing}

For $\alpha =\frac{7}{2}$ we will now present the plots of the waves
for an increasing sequence of $\beta$ values i.e. 4, 8, 12, and
20 (in order to investigate the ``infinite beta limit'' already
introduced in our previous work \cite{7}). We will compute the function
\begin{equation}
\psi ( k;\frac{7}{2},\beta ,\frac{1}{2}) =k^{\frac{3}{\beta }}c_{k}(
\frac{7}{2},\beta ) 
\end{equation}

\noindent which will also be compared with the expression given
by the Baez-Duarte formula (19) in the asymptotic region $k\overset{
}{\rightarrow }\infty $. Here we will take into account only the
contribution of the groundstate of the spectrum i.e $z=\frac{1}{2}+i
14.134725141$. It is then convenient to introduce the new variable
$x=\log ( k) $. This allow us to control more efficiently the wavelength
and the amplitude of the wave in the region to be considered ({\itshape
x} runs from 8 to 22, so {\itshape k} up to $3.6\cdot {10}^{9}$).

In the Figures 8-11 we present our numerical results for increasing
$\beta$ values, which we call the ``strong coupling limit''.
\begin{figure}[h]
\begin{center}
\includegraphics{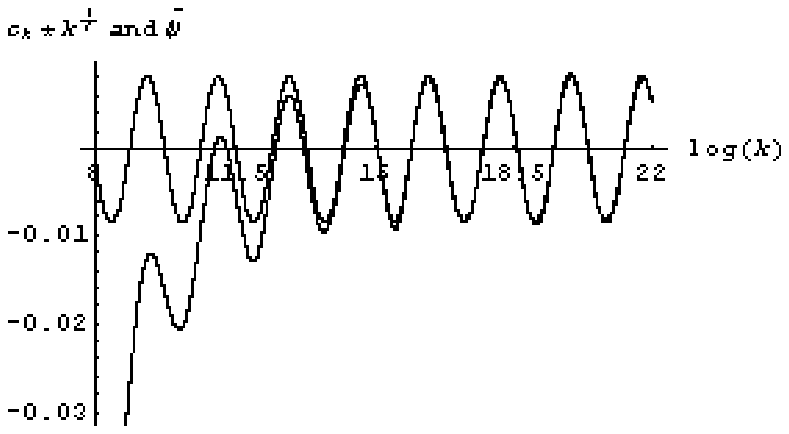}

\end{center}
\caption{The wave $k^{\frac{3}{4}}c_{k}$ (lowest curve) and the
approximation $\overline{\psi }$ (highest curve), $\beta =4$}

\end{figure}
\begin{figure}[h]
\begin{center}
\includegraphics{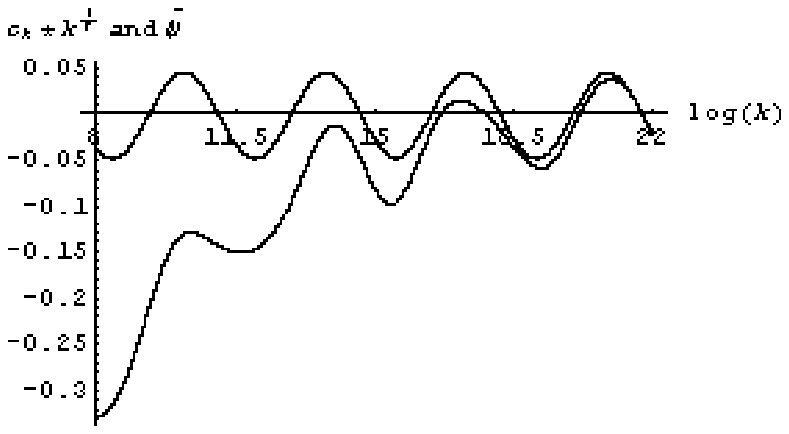}

\end{center}
\caption{The wave $k^{\frac{3}{8}}c_{k}$ (lowest curve) and the
approximation $\overline{\psi }$ (highest curve), $\beta =8$}

\end{figure}
\begin{figure}[h]
\begin{center}
\includegraphics{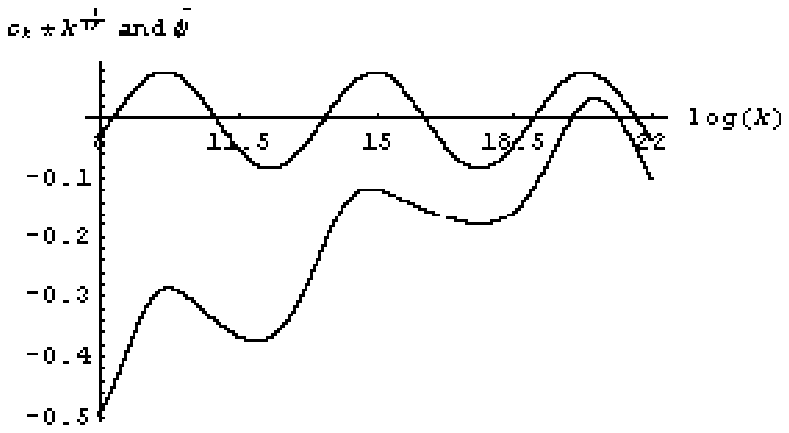}

\end{center}
\caption{The wave $k^{\frac{3}{12}}c_{k}$ (lowest curve) and the
approximation $\overline{\psi }$ (highest curve), $\beta =12$}

\end{figure}
\begin{figure}[h]
\begin{center}
\includegraphics{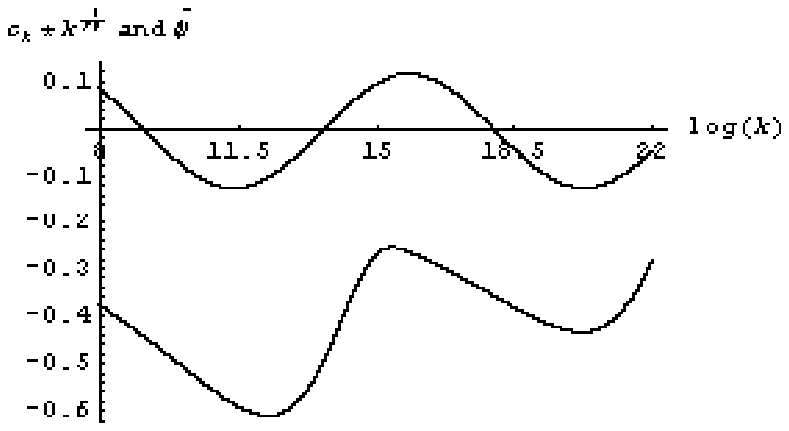}

\end{center}
\caption{The wave $k^{\frac{3}{20}}c_{k}$ (lowest curve) and the
approximation $\overline{\psi }$ (highest curve), $\beta =20$}

\end{figure}

At the same time it is seen that in this case $|c_{k}|$ itself is
smaller than ($c_{k}$ is not the critical function!):
\begin{equation}
\left| \frac{1}{\zeta ( \frac{7}{2}) }-1\right| \cong 0.11247897
\end{equation}

\noindent at least for the case $\beta =4$ as already discussed
in our previous work \cite{7} concerning only very low values of
{\itshape k}. Figure 12 confirm this behaviour also for large value
of {\itshape k}. For this example the region of annihilation of
the ``eincoming`` wave extends up to larger and larger values of
{\itshape k}. It should be noted that for the critical function
$\psi$ (20) the situation is more delicate since the value of a
possible bound on $\psi$ depends on $\beta$ as it is been from our
numerical results.
\begin{figure}[h]
\begin{center}
\includegraphics{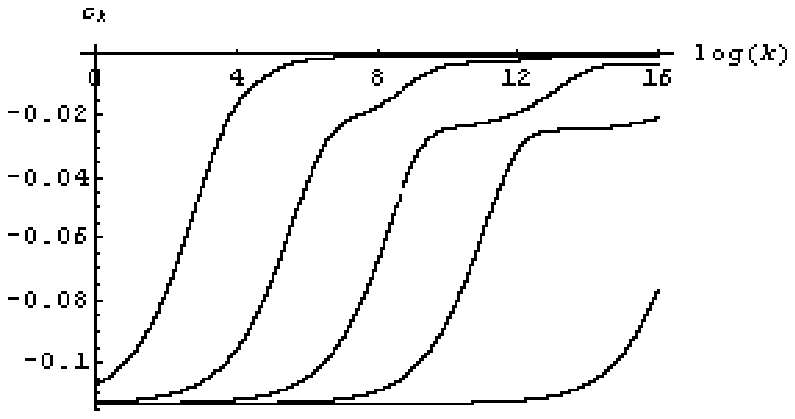}

\end{center}
\caption{$c_{k}( \frac{7}{2},\beta ) $ for $\beta =4,8,12,16,24$
(from left to right)}

\end{figure}

\section{Conclusions}

In this work we have extended the numerical experiments obtained
in our previous work and considered many cases of ``waves'' (with
two parameters $\alpha$ and $\beta$). The experiments allow to obtain
up to 6 oscillations at low $\beta$ values for (20) whose amplitudes
are in agreement with those given by an extension of an asymptotic
formula due to Baez-Duarte for the case $\alpha =\beta =2$.

In the process of increasing the values of $\beta$ (at least in
the case $\alpha =\frac{7}{2}$) it has been shown that the width
of the ``annihilation'' region increase with $\beta$ where the amplitude
of the wave seems (still) to remain bounded. In this connection
we may argue that the results of our numerical experiments indicate
that RH may barely be true due to the behaviour of $\psi$ in the
large $\beta$ limit. It is also conjectured that an absolute bound
on $\psi$, for all $\alpha$ and $\beta$ allowed, should be given
by $|\frac{1}{\zeta ( \frac{1}{2}+\epsilon ) }-1|\cong 1.68$. In
fact Littlewood has shown that on RH, $\sum \limits_{n=1}^{\infty
}\frac{\mu ( n) }{n^{\frac{1}{2}+\epsilon }}$ is convergent for
all $\epsilon >0$. In the light of the results of the experiments
obtained so far and with the understanding that our remark is speculative,
we belive that even with more sofisticated experiments it will be
very hard to obtain values of the critical function $k^{a}c_{k}$
which in absolute value will be greater at large {\itshape k} then
those of the infinite $\beta$ limit, as commented in the Appendix.

\appendix

\section{}

In the context of the numerical experiments performed so far, it
is helpful to obtain a crude inequality concerning a bound on the
critical function. This is simply obtained by setting\ \ $|\mu (
n) |=1$ in (3). The critical function in the representation of $\frac{1}{\zeta
( s) }$ in terms of the two parameter Pochammer's polynomials is
given by:
\[
k^{\frac{\alpha -\frac{1}{2}}{\beta }}c_{k}\cong k^{\frac{\alpha
-\frac{1}{2}}{\beta }}\sum \limits_{n=1}^{N}\frac{\mu ( n) }{n^{\alpha
}}{\left( 1-\frac{1}{n^{\beta }}\right) }^{k}=:f_{k}( \alpha ,\beta
,N) 
\]

\noindent where {\itshape N} is the maximum value of the argument
in the M\"obius function considered in a ideal numerical experiment
({\itshape N} finite). We then have, introducing the variable $x=\log
( k) $ that:
\[
\left| f_{k}( \alpha ,\beta ,N) \right| \leq e^{\frac{\alpha -\frac{1}{2}}{\beta
}x}( \zeta ( \alpha ) -1) e^{\log ( 1-\frac{1}{N^{\beta }}) e^{x}}
\]

For large {\itshape N} we have:
\[
\left| f_{k}( \alpha ,\beta ,N) \right| \leq \left( \zeta ( \alpha
) -1\right) e^{\frac{\alpha -\frac{1}{2}}{\beta }x-\frac{1}{N^{\beta
}}e^{x}}
\]

As an example we consider our case $\alpha =\frac{7}{2}$ and $\beta
=4$. Remembering that from Table 2 the amplitude calculated only
with the first non trivial zero is about 0.008411, we may ask: for
what {\itshape N} and {\itshape k,} $|f_{k}( \alpha ,\beta ,N) |$
is bounded by the value 0.008411? For example the inequality is
satisfied for the followig pairs:
\[
\begin{array}{cc}
 N=1000 & x>31 \\
 N={10}^{6} & x>60
\end{array}
\]

As a second example we consider the Riesz case ($\alpha =\beta =2$).
From the Table 2, the amplitude (still restricting to the contribution
of the first zero) is 0.000078. The inequality is satisfied as follows:
\[
\begin{array}{cc}
 N=1000 & x>17 \\
 N={10}^{6} & x>31 \\
 N={10}^{9} & x>87.2
\end{array}
\]

Now, still for the Riesz case
\[
\left| \frac{1}{\zeta ( s) }\right| \leq \sum \limits_{k=0}^{\infty
}\left| P_{k}( s,2,2) c_{k}( 2,2) \right| 
\]

\noindent assuming as seen in the experiments that $k^{\frac{3}{2\beta
}}c_{k}( 2,\beta ) $, i.e. the critical function is bounded then
we would have that:
\[
\left| \frac{1}{\zeta ( s) }\right| \leq \sum \limits_{k=0}^{\infty
}\frac{1}{k^{1+\frac{\delta }{2}}}k^{\frac{3}{2\beta }}c_{k}( 2,\beta
) \leq C
\]

\noindent Thus, $\frac{1}{\zeta ( s) }$ would be different from
zero for $\mathfrak{R}( s) >\frac{1}{2}+\delta  , \delta >0$. It
should be said that for $\alpha =\frac{3}{2}$ and $\beta =1$, in
the numerical experiments, the wave seems to decay at zero after
a few of oscillations, a situation very different from the case
$\alpha =\frac{7}{2}$ and $\beta =4$ (but the above inequality apply
also in this case).

So, if the critical function for any $\beta$ is in absolute values
bounded for $k>K$ by the value of the infinite $\beta$ limit 1.68\ldots
, the RH should be true.

In the plot of $\psi$ for the case $\alpha =\frac{7}{2}$ and $\beta
=4$ (Figure 13) up to $k=10$ billions we see 8 oscillations.
\begin{figure}[h]
\begin{center}
\includegraphics{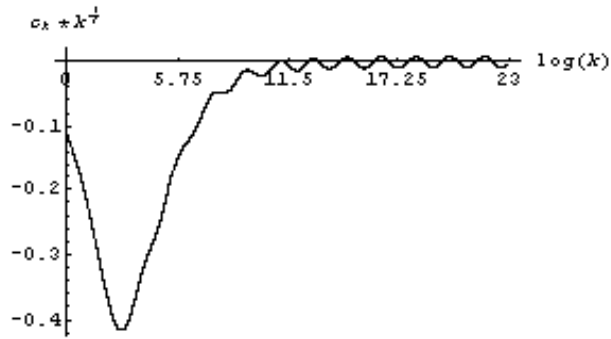}

\end{center}
\caption{$k^{\frac{3}{4}}c_{k}$ for $\alpha =\frac{7}{2}$ and $\beta
=4$ until $k=10\cdot {10}^{9}$\ \ }

\end{figure}

Final comments.
\begin{enumerate}
\item Further\ \ application of the crude inequality considered
in this Appendix indicates that in a ideal experiment using (3),
with $N={10}^{10}$, the amplitude of the wave for $\log ( k) <23$
will change at most ${10}^{-6}$ time the value 0.008411... obtained
with $N={10}^{6}$ in (3) and compatible with the Baez-Duarte amplitude
using only the first nontrivial zero of Zeta. This indicates some
stability of the numerical experiments in the intermediate range
$\log ( k) <23$ (see Figure 13).
\item In the general case of an ideal experiment with large {\itshape
N} in (3) the crude inequality indicates also that for $k>K( N)
$, {\itshape K} finite, the amplitude of the critical function $\psi$
is bounded in absolute value by $|\frac{1}{\zeta ( \frac{1}{2})
}-1|\cong 1.68$, choosing {\itshape N} as needed. In the same way
we may argue that the amplitude of the critical function may be
obtained as small as we want (in particular smaller then 0.008411...)
choosing {\itshape N} as needed for values of $k>K( N) $, {\itshape
K} still finite.
\item One of the open questions is now the following: the critical
function at large value of {\itshape k} is growing, stabilizing
to a ``periodic pure wave'' with constant amplitude or decaying
with a zero amplitude? From the results of our numerical treatment
we are more in favour of the last two cases.\ \ \ 
\end{enumerate}

\end{document}